\documentclass[journal, twoside]{IEEEtran}

\usepackage{import}

\usepackage[utf8]{luainputenc}

\usepackage[T1]{fontenc}

\usepackage{graphicx}
\usepackage[dvipsnames]{xcolor}

\usepackage[english]{babel}
\addto\captionsenglish{}
\addto\captionsenglish{}
\usepackage{csquotes}

\usepackage{newtxtext}
\usepackage{amsthm}
\usepackage[slantedGreek]{newtxmath}
\usepackage[OMLmathsfit]{isomath}
\DeclareMathAlphabet{\mathbfsf}{\encodingdefault}{\sfdefault}{bx}{n}
\usepackage{bm}
\usepackage{envmath}
\usepackage{mathtools}
\usepackage{commath}
\usepackage{siunitx}
\DeclareSIUnit\ohm{\upOmega}
\usepackage{yfonts}

\usepackage{caption}
\usepackage{subcaption}
\graphicspath{{figures/}} 

\usepackage{booktabs}
\usepackage{footmisc}  
\usepackage{multirow}
\usepackage{makecell}

\usepackage{algpseudocode}
\usepackage{algorithm}

\usepackage[style=ieee, backend=biber, isbn=false, maxbibnames=9]{biblatex}

\usepackage{comment}

\usepackage{url}
\usepackage{nicefrac}

\theoremstyle{definition}

\theoremstyle{plain}

\theoremstyle{remark}

\usepackage{lineno}
\modulolinenumbers[5]
\usepackage{todonotes}
\usepackage{umoline}

\usepackage{pgfplots}
\usepackage{pgfplotstable}
\pgfplotsset{compat=newest}
\pgfplotsset{plot coordinates/math parser=false}
\newlength\figureheight
\newlength\figurewidth
\pgfplotsset{every axis plot/.append style={line width=1.5pt},
    legend style={font=\footnotesize, 
        text height=1.0ex,
        draw=black,
        fill=white,
        legend cell align=left}}
\usetikzlibrary{arrows.meta}
\usepackage[hidelinks]{hyperref} 
\usepackage[english]{cleveref}

\Crefname{defn}{definition}{definitions}
\Crefname{defn}{Definition}{Definitions}

\Crefname{asm}{assumption}{assumptions}
\Crefname{asm}{Assumption}{Assumptions}

\crefname{lem}{lemma}{lemmas} 
\Crefname{lem}{Lemma}{Lemmas}

\crefname{prop}{proposition}{propositions} 
\Crefname{prop}{Proposition}{Propositions}

\crefname{thm}{theorem}{theorms} 
\Crefname{thm}{Theorem}{Theorms}

\crefname{cor}{corollary}{corollaries}
\Crefname{cor}{Corollary}{Corollaries}

\crefname{figure}{Fig.}{Figs.}
\Crefname{figure}{Fig.}{Figs.}
\newcounter{subequation}
\newlength\mtabskip\mtabskip=-1.25cm

\def\mtabLong{long}

\newcommand{\mr}{\mathrm}

\newcommand{\mc}{\mathcal}

\newcommand{\veg}[1]{\bm{#1}}     
\newcommand{\mat}[1]{\mathsfbfit{#1}} 
\renewcommand{\vec}[1]{\mathsfbfit{#1}} 
\newcommand{\wvec}[1]{\widetilde{\mathsfbfit{#1}}} 
\newcommand{\vecop}[1]{\bm{\mathcal{#1}}} 

\newcommand{\dyd}[1]{\bm{#1}}   

\newcommand{\n}{\hat{\bm{n}}}

\renewcommand{\Re}{\mathrm{Re}}
\renewcommand{\Im}{\mathrm{Im}}

\newcommand{\dd}{\mathrm{d}}  

\newcommand{\jm}{\mathrm{j}}  

\newcommand{\e}{\mathrm{e}}

\newcommand{\R}{\mathbb{R}}

\DeclareMathOperator{\nablat}{\nabla_\mathrm{t}}

\newcommand{\T}{\mr{T}}

\newcommand{\bJ}{\operatorname{J}}
\newcommand{\bY}{\operatorname{Y}}
\newcommand{\bH}{\operatorname{H}}

\newcommand\restr[2]{{
        \left.\kern-\nulldelimiterspace 
        #1 
        \vphantom{|} 
        \right|_{#2} 
}}

\newcommand\rst[3]{{
        \left.\kern-\nulldelimiterspace 
        #1 
        \vphantom{|} 
        \right|_{#2}^{#3} 
}}

\newcommand{\gf}{g^\mathrm{f}}
\newcommand{\gaux}{g^\mathrm{aux}}
\newcommand{\gfH}{\veg{g}_H^\mr{f}}
\newcommand{\gauxH}{\veg{g}_H^\mr{aux}}
 
\usepackage{acro}

\DeclareAcronym{DG}
{
    short = DG ,
    long = discontinuous Galerkin
}

\DeclareAcronym{ACA}
{
    short = ACA ,
    long = adaptive cross approximation
}

\DeclareAcronym{EFIE}
{
    short =  EFIE ,
    long = electric field integral equation
}

\DeclareAcronym{MFIE}
{
    short =  MFIE ,
    long = magnetic field integral equation
}

\DeclareAcronym{MUIE}
{
    short =  MUIE ,
    long = Müller integral equation
}

\DeclareAcronym{PMCHWT}
{
    short =  PMCHWT ,
    long = Poggio-Miller-Chang-Harrington-Wu-Tsai integral equation
}

\DeclareAcronym{SPD}
{
    short =  SPD ,
    long = {symmetric, positive definite}
}

\DeclareAcronym{SPSD}
{
    short =  SPD ,
    long = {symmetric, positive semi-definite}
}

\DeclareAcronym{PEC}
{
    short =  PEC ,
    long = perfectly electrically conducting
}

\DeclareAcronym{PMC}
{
    short =  PMC ,
    long = perfectly magnetically conducting
}

\DeclareAcronym{RWG}
{
    short = RWG ,
    long = Rao-Wilton-Glisson
} 

\DeclareAcronym{BC}
{
    short = BC ,
    long = Buffa-Christiansen
}

\DeclareAcronym{SVD}
{
    short = SVD ,
    long = singular value decomposition
}

\DeclareAcronym{CG}
{
    short = CG ,
    long = conjugate gradient
} 

\DeclareAcronym{PCG}
{
    short = PCG ,
    long = preconditioned conjugate gradient
} 

\DeclareAcronym{CGS}
{
    short = CGS ,
    long = conjugate gradient squared
}

\DeclareAcronym{CMP}
{
    short = CMP ,
    long = Calderón multiplicative preconditioner
} 

\DeclareAcronym{RFCMP}
{
    short = RF-CMP ,
    long = refinement-free Calderón multiplicative preconditioner
} 

\DeclareAcronym{HPD}
{
    short = HPD ,
    long = {Hermitian, positive definite}
} 

\DeclareAcronym{RHS}
{
    short = RHS ,
    long = right-hand side
}

\DeclareAcronym{LSE}
{
    short = LSE ,
    long = linear system of equations
}

\DeclareAcronym{AMG}
{
    short = AMG ,
    long = algebraic multigrid
}

\DeclareAcronym{PW}
{
    short = PW ,
    long = plane wave
}

\DeclareAcronym{GMRES}
{
    short = GMRES ,
    long = generalized minimum residual
}

\DeclareAcronym{IDR}
{
    short = IDR ,
    long = induced dimension reduction
}

\DeclareAcronym{BICGstab}
{
    short = BiCGstab ,
    long = stabilized bi-conjugate gradient
}

\DeclareAcronym{FF}
{
    short = FF ,
    long = far field
}

\DeclareAcronym{NF}
{
    short = NF ,
    long = near field
}

\DeclareAcronym{MoM}
{
    short = MoM ,
    long = method of moments
}

\DeclareAcronym{SIE}
{
    short = SIE ,
    long = surface integral equation
}

\DeclareAcronym{GSIE}
{
    short = GSIE ,
    long = generalized source integral equation
}

\DeclareAcronym{LR}
{
    short = LR ,
    long = low-rank
}

\DeclareAcronym{MGD}
{
    short = MGD ,
    long = modified Green’s dyadic
}

\DeclareAcronym{GEIE}
{
    short = GEIE ,
    long = generalized equivalence IE
}

\DeclareAcronym{CSIE}
{
    short = CSIE ,
    long = combined source IE
}

\DeclareAcronym{CFIE}
{
    short = CFIE ,
    long = combined field IE
}

\DeclareAcronym{TMGSIE}
{
    short = TM-GSIE ,
    long = TM \ac{GSIE}
}

\DeclareAcronym{TE}
{
    short =  TE ,
    long = transverse electric
}

\DeclareAcronym{TM}
{
    short =  TM ,
    long = transverse magnetic
}

\DeclareAcronym{MGF}
{
    short =  MGF ,
    long = modified Green's function
}

\DeclareAcronym{EV}
{
    short = EV,
    long = eigenvalue
}

\DeclareAcronym{Y}
{
    short = Y,
    long = Yukawa
}

\DeclareAcronym{GC}
{
    short = GC,
    long = generalized Calderón
} 

\newcolumntype {n}{c}
\newcolumntype {N}{>{\small}c}
\newcolumntype {L}{>{\small}l}
\newcolumntype {F}{>{\footnotesize}c}
\newcolumntype {v}[1]{>{\raggedright \hspace {0pt}} p {#1}}
\newcolumntype {V}[1]{>{\small \raggedright \hspace {0pt}} p {#1}}
\newcolumntype{d}[1]{>{\DC@{.}{.}{#1}}c<{\DC@end}}

\newcolumntype{R}[1]{%
    >{\begin{turn}{90}\begin{minipage}{#1}\small\raggedright\hspace{0pt}}l%
            <{\end{minipage}\end{turn}}%
}

\NewDocumentCommand{\TA}{o}{
    \IfNoValueTF {#1} {%
        \vecop T_{\kern-2pt\mr{A}}
    }
    {
        \vecop T_{\kern-2pt\mr{A},#1}
    }
}

\NewDocumentCommand{\TPhi}{o}{
    \IfNoValueTF {#1} {%
        \vecop T_{\kern-2pt\Phiup}
    }
    {
        \vecop T_{\kern-2pt\Phiup,#1}
    }
}

\NewDocumentCommand{\matTA}{o}{
    \IfNoValueTF {#1} {%
        \mat T_\mr{A}   
        }
    {
        \mat T_{\mr{A},#1}
    }
}

\NewDocumentCommand{\matTPhi}{o}{
    \IfNoValueTF {#1} {%
        \mat T_\Phiup   
        }
    {
        \mat T_{\Phiup,#1}
    }
}

\NewDocumentCommand{\MSL}{o}{
    \IfNoValueTF {#1} {%
        \veg \Psi_\mr{SL}
        }
    {
        \veg \Psi_{\mr{SL},#1}
    }
}

\NewDocumentCommand{\MDL}{o}{
    \IfNoValueTF {#1} {%
        \veg \Psi_\mr{DL}
        }
    {
        \veg \Psi_{\mr{DL},#1}
    }
}

\NewDocumentCommand{\PA}{o}{
    \IfNoValueTF {#1} {%
        \veg \Psi_\mr{A}
        }
    {
        \veg \Psi_{\mr{A},#1}
    }
}

\NewDocumentCommand{\PPhi}{o}{
    \IfNoValueTF {#1} {%
        \veg \Psi_{\Phiup}
        }
    {
        \veg \Psi_{\Phiup,#1}
    }
}

\NewDocumentCommand{\TM}{o}{
    \IfNoValueTF {#1} {%
        \vecop T_\mathrm{TM}
        }
    {
         \vecop T_\mathrm{TM, #1}
    }
}

\NewDocumentCommand{\TMa}{o}{
    \IfNoValueTF {#1} {%
        \vecop T_\mathrm{TM}^{\kern+2pt\mathrm{aux}}
        }
    {
         \vecop T_\mathrm{TM, #1}^{\kern+2pt\mathrm{aux}}
    }
}

\NewDocumentCommand{\TMm}{o}{
    \IfNoValueTF {#1} {%
        \vecop T_\mathrm{TM}^{\kern+2pt\mathrm{m}}
        }
    {
         \vecop T_\mathrm{TM,#1}^{\kern+2pt\mathrm{m}}
    }
}

\NewDocumentCommand{\TE}{o}{
    \IfNoValueTF {#1} {%
        \vecop T_\mathrm{TE}
        }
    {
         \vecop T_\mathrm{TE,#1}
    }
}

\NewDocumentCommand{\TEa}{o}{
    \IfNoValueTF {#1} {%
        \vecop T_\mathrm{TE}^{\kern+2pt\mathrm{aux}}
        }
    {
         \vecop T_\mathrm{TE,#1}^{\kern+2pt\mathrm{aux}}
    }
}

\NewDocumentCommand{\TEm}{o}{
    \IfNoValueTF {#1} {%
        \vecop T_\mathrm{TE}^{\kern+2pt\mathrm{m}}
        }
    {
         \vecop T_\mathrm{TE,#1}^{\kern+2pt\mathrm{m}}
    }
}

\newcommand{\hveg}[1]{\hat{\veg{#1}}}
\newcommand{\wveg}[1]{\widetilde{\veg{#1}}}

\renewcommand{\wvec}[1]{\widetilde{\vec{#1}}}

\begin{document}

	\title{Well-conditioned Electric Field Surface Integral Equations using Reflective Generalized Sources}

	%
	
	\author{Yossi Dahan,~\IEEEmembership{Student Member,~IEEE,}
	        Suryakumar Sivakumar,~\IEEEmembership{Student Member,~IEEE,}\\
	        Yaniv Brick,~\IEEEmembership{Senior Member,~IEEE,}
			and~Simon B.~Adrian,~\IEEEmembership{Senior Member,~IEEE}
	\thanks{Funded by the Deutsche Forschungsgemeinschaft (DFG, German Research Foundation)~--~544296727. \emph{(Yossi Dahan and Suryakumar Sivakumar are co-first authors.)} \emph{(Corresponding author: Yaniv Brick.)}}%
    \thanks{Y.~Dahan and Y.~Brick are with the School of Electrical and Computer Engineering, Ben-Gurion University of the Negev, 8410501, Beer-Sheva, Israel (e-mail: ybrick@bgu.ac.il).}
	\thanks{S.~Sivakumar and S.~B.~Adrian are with the Fakultät für Informatik und Elektrotechnik, Universität Rostock, 18059 Rostock, Germany (e-mail: simon.adrian@uni-rostock.de).}
	}

	\maketitle

	\begin{abstract}
        This work uses the generalized source approach to develop a class of well-conditioned integral operators that are free of internal resonance, without the need for combined formulations. 
        The \Acp{GSIE} kernels are obtained by augmenting the conventional \ac{EFIE} kernel with auxiliary contributions to enhance the rank deficiency of the corresponding moment matrix blocks.
        This paper presents the first investigation of the spectra of \ac{GSIE} operators for auxiliary kernels produced by internal scattering convex shields. 
        Using closed-form expressions for concentric circular scatterers and shields, it is shown that, with shield parameters that are suitable for enhanced compressibility, the \ac{TM}- and \ac{TE}-\ac{GSIE} operators are free of internal proper and quasi-resonances.
        The auxiliary components are shown to be compact perturbations of their \ac{EFIE} counterparts. 
        Hence these \acp{GSIE} inherit their dense-discretization breakdown, which remains curable via Calderón-type preconditioning.
        The mechanisms that govern the high-frequency breakdown are shown to be influenced by the auxiliary component, leading, in some cases, to greater resilience. 
        These observations are shown to remain valid for moment matrices and \acp{GSIE} designed with non-circular stencil shields.
        For both \ac{GSIE}-dual and Yukawa-kernel preconditioners, the formulations exhibit favorable spectral properties while maintaining their compressibility.
        This makes the formulations attractive for the design of fast iterative solvers.
        The results on the two-dimensional shield-based operators provide a foundation for extending the approach to three-dimensional problems and to auxiliary kernels with broader geometric applicability.
	\end{abstract}
\acresetall


	%
	\IEEEpeerreviewmaketitle
    
    \section{Introduction}
    \IEEEPARstart{S}{URFACE} integral equations (\acsp{SIE}) are commonly used for the scattering analysis of impenetrable objects~\cite{jin_theory_2015}.
    In their discretization, for example by using the \ac{MoM}~\cite{harrington_field_1993}, they give rise to systems of linear equations described by dense matrices.
    The ability to accurately solve these systems within reasonable computation times and resources relies on two components: i) an efficient representation of the system matrix that enables its fast matrix--vector products (in the context of iterative solvers 
    \cite{chew_fast_2001, bleszynski1996aim, wei2013more, brick2009multilevel, michielssen1996multilevel, zhao2005adaptive, tamayo2011multilevel}) and factorization (as part of direct solvers \cite{martinsson2005fast, martinsson2007fast, shaeffer2008direct, brick2011fast, chai2012direct, wei2012fast, corona2015n, guo2017butterfly, heldring2025fast})
    to avoid poor scaling of the computational costs with the number of unknowns $N$, and ii) a well-conditioned discretized \acs{SIE} formulation (see \cite{adrian_ElectromagneticIntegralEquations_2021} and references therein), which ensures both a low algebraic solution error and a small number of iterations $N_{\text{iter}}$ for solver convergence.  
    These are especially critical for multiscale problems, where the analyzed structures are both electrically large and include small and finely discretized features and, therefore, involve mixtures of ill-conditioning mechanisms.
    With large $N$, fast iterative solvers have been favored for the solution of such problems, aiming to achieve $N_{\text{iter}}\mathcal{O}(N\log N)$ computational costs with a size-independent $N_{\text{iter}}$.
         
    Fast solvers can be classified as methods that rely on either “matrix-free” representations of field integrals (e.g., \cite{chew_fast_2001, bleszynski1996aim, wei2013more, brick2009multilevel}) or algebraic compression of admissible moment matrix blocks (e.g., \cite{zhao2005adaptive, tamayo2011multilevel, michielssen1996multilevel, brick2011fast, hackbusch1999sparse, borm2007data, tetzner_IncompleteAdaptiveCross_2026}). 
    The main appeal of the latter class lies in their kernel-independent implementation as well as their role as a foundation for fast direct solvers, when accompanied by a suitable compressed matrix arithmetic \cite{hackbusch1999sparse, borm2007data, liu2021butterfly, heldring2025fast}.
    Within this class, \ac{LR}-compression-based methods are particularly convenient \cite{martinsson2005fast, martinsson2007fast, shaeffer2008direct, brick2011fast, chai2012direct, wei2012fast, corona2015n, zhao2005adaptive, tetzner_IncompleteAdaptiveCross_2026}. 
    Compression and speedup were demonstrated by these techniques for geometrically complex and electrically large scatterers. 
    However, only for objects where the interactions are of reduced geometric dimensionality (e.g., in quasi-planar \cite{martinsson2005fast, martinsson2007fast, corona2015n} and elongated structures \cite{ martinsson2007fast, corona2015n}), they also exhibit slow asymptotic scaling of block ranks and associated computational costs;
    for arbitrary geometries, the ranks of interactions between large subdomains, via the kernels of conventional \acp{SIE} (e.g., \ac{EFIE} and \ac{MFIE}), are governed by their broadside components (i.e., broadside with respect to the surface subdomain of the current densities)~\cite{brick2026interpreting, brick2016fast, brick2018rapid}.
    As a result, their asymptotic scaling remains linear in $N$ and savings are limited to constant (often large) factors.

    Conventional \acp{SIE}, specifically the \ac{EFIE}, are also plagued by mechanisms that cause deterioration of the condition number (or ill-conditioning) \cite{saad_iterative_2003}.
    These include interior resonances, which are typically eliminated by combining \acp{SIE}, as well as the dense-discretization \cite{andriulli_MultiplicativeCalderonPreconditioner_2008} and high-frequency \cite{darbas_PreconditioneeursAnalytiquesType_2004, boubendir_WellconditionedBoundaryIntegral_2014} breakdowns, which are commonly combated by using preconditioners \cite{christiansen_PreconditionerElectricField_2002, andriulli_MultiplicativeCalderonPreconditioner_2008,vipiana_EFIEModelingHighdefinition_2010,andriulli_WellconditionedElectricField_2013,boubendir_WellconditionedBoundaryIntegral_2014, carpentieri_PreconditioningLargescaleBoundary_2014, adrian_HierarchicalPreconditionerElectric_2017, adrian_RefinementfreeCalderonMultiplicative_2019, adrian_ElectromagneticIntegralEquations_2021,antoine_IntroductionOperatorPreconditioning_2021,fierro-piccardo_OSRCPreconditionerEFIE_2023,darbas_GeneralFrameworkOSRCpreconditioned_2025}.
    Fast solution of multiscale problems through \ac{LR}-compression-based methods requires formulations that are resonance free, immune to the dense-discretization breakdown, well-conditioned for broad frequency ranges, and exhibit slow scaling of interaction ranks with increasing electrical size.
    
    For essentially convex geometries, \ac{LR}-compressibility can be significantly improved. 
    This is achieved by modifying the \ac{SIE} kernel to attenuate the broadside interactions between opposing surfaces.
    This accentuates the endfire component of the interactions, which appear to have a reduced effective geometric dimensionality~\cite{boag2012generalized, brick2014fastEssentially, sharshevsky2020direct, zvulun2023generalized, dahan2024fast, dahan2025vector, kalhofer2025fast, kalhofer2026multipole}, leading to slower scaling of the ranks with $N$.
    \Acp{GSIE} replace the free-space Green’s function kernel with \acp{MGF} that incorporate auxiliary components into the conventional \ac{EFIE} kernel.
    These contributions can be attributed to auxiliary sources that are internal to the scatterer. 
    They are designed to approximately cancel the radiation of the original source into the scatterer volume and toward the opposite side.
    To this end, various mechanisms have been proposed, including reflective (see~\cite{boag2012generalized, brick2014fastEssentially, dahan2024fast, dahan2025vector})
    and absorptive (see~\cite{sharshevsky2020direct, zvulun2023generalized}) “shield” surfaces  
    and auxiliary multipole sources~\cite{kalhofer2025fast, kalhofer2026multipole}.
    Several kernel designs, differing in the width and depth of the attenuation region, have been successfully used to improve compressibility in various fast direct solver architectures.
    However, well-conditioned \ac{GSIE} formulations suitable for fast iterative solution of multiscale problems are yet to be presented. 
    
    This work presents a class of well-conditioned \acp{GSIE} for 2-D scattering problems.
    Their development builds upon a study of the spectral properties of existing \acp{GSIE} formulations.
    As \ac{GSIE} operators are formulated as perturbations of \ac{EFIE} operators, they are expected to inherit many spectral properties of their conventional counterparts.   
    Reflective-type \acp{GSIE} provide a convenient setting for examining this assumption.
    In these formulations, the auxiliary kernel components can be attributed to the reflection of the free-space \ac{SIE} kernel field by smooth, convex, impenetrable surfaces internal to the scatterer~\cite{brick2014fastEssentially,dahan2024fast}.
    For the particular class of circular shields~\cite{brick2014fastEssentially, brick2014fast}, \ac{MGF} can be written analytically as trigonometric eigenfunction expansions. 
    Moreover, for configurations of concentric circular scatterer and coinciding shields, the operator spectra admit closed-form expressions.
    These spectra are obtained here using techniques known for standard electromagnetic operators (see \cite{hsiao_ErrorAnalysisNumerical_1994, vico_BoundaryIntegralEquation_2014} and references therein).
    Examining the \ac{TM} \cite{dahan2024fast} and \ac{TE} \cite{dahan2025vector} \ac{GSIE} operator spectra, suitable design parameters, for which the \acp{GSIE} are practically resonance-free, are identified.
    The closed-form spectra further reveal that the scaling with mesh size $h$ of the condition numbers of these operators is governed by that of the conventional kernels. 
    Accordingly, both the analysis and the numerical results show that the inherited dense-discretization breakdown can be cured by using established Calderón-type preconditioners \cite{cools2013spectrum}.
    The numerical results also show that these properties extend to non-circular scatterers and stencil shields.
    
    The remainder of this paper is organized as follows: \Cref{sec:background} reviews the background required for this work, including the \ac{TM} and \ac{TE} electric field \ac{GSIE} formulations using reflective shields.
    \Cref{sec:Analysis} analyzes the \acp{GSIE} spectra and conditioning for the special case of a concentric circular scatterer and shield.
    \Cref{sec:NumResults} compares these predictions with numerical results obtained from \ac{MoM} discretizations.
    \Cref{sec:Conclusion} summarizes the principal findings.
    
    \section{Background}\label{sec:background}
    \subsection{SIE Formulation of 2-D Scattering Problems}
    
    The scattering of a time-harmonic field $\veg e^\text{inc}$ (an $\e^{\jm{\omega}t}$ dependence is assumed and suppressed) by a \ac{PEC} surface $S$, with normal and tangent vectors $\n$ and $\hveg{t}$ (as shown in \Cref{fig:scatteringscenario}), can be described using \acp{SIE}.
    For instance, in the conventional \ac{EFIE}, the scattered field is described by an integral over the induced physical surface current $\veg j$, that is,
    \begin{equation}
    \hveg{n}\times\[-\jm k\int_S \dyd G^\mr{f}(\veg r, \veg r')\cdot\veg j(\veg r')\, \dd t'\] = -\frac{1}{\eta}\n \times \veg e^\text{inc}\, , \quad \veg r\in S\, . \label{eq:EFIE}
    \end{equation}
    Here, $\dyd G^\mr{f}(\veg r,\veg r') \coloneq (\dyd I+1/k^2\,\nabla\nabla) g^\mr{f}(\veg r, \veg r') $ and $g^\mr{f}(\veg r,\veg r')$ are the free-space Green's dyadic and scalar Green’s function, $\dyd I$ is the identity dyadic, $\eta$ is the free-space impedance, and $k=2\uppi/\lambda$ is the wave number.
    For two-dimensional (2-D) $z$-invariant \ac{TM} problems, the \ac{EFIE} in \eqref{eq:EFIE} reduces to a scalar equation for the $z$-components
    \begin{equation}
    \TM(j_\mathit{z})(\veg r) =-\jm k\int_S g^\mr{f}(\veg r,\veg r') j_\mathit{z}(\veg r')\, \dd t'=-\frac{1}{\eta}e_\mathit{z}^\text{inc}(\veg r)\, , \label{eq:TM-EFIE}
    \end{equation}
    where $g^\mr{f}(\veg r,\veg r')=-\jm\bH_0^{(2)}(k|\veg r-\veg r'|)/4$ with $\bH_0^{(2)}(\cdot)$ the zeroth order Hankel function of the second kind.
    For $z$-invariant \ac{TE} problems, the Green’s dyadic reduces to
    \begin{equation}
    \dyd G^\mr{f}(\veg r,\veg r')=\(\dyd I+\frac{1}{k^2}\nablat\nablat\)g^\mr{f}(\veg r,\veg r')\, , \label{eq:TE-GD}
    \end{equation}
    where dyadic notation indicates $2\times2$ entities and $\nablat$ is the transverse gradient.
    Let $j_\mr{t}(\veg r)$ and $e_\mr{t}^\text{inc}(\veg r)$ denote the surface current and the tangential component of the incident transverse electric field on $S$, and let $\uppartial_\mr{t}$ denote the tangential derivative.
    The TE-\ac{EFIE} operator can be written as
    \begin{multline}
    \TE(j_\mr{t})(\veg r)=-\jm k\int_S g^\mr{f}(\veg r,\veg r') j_\mr{t}(\veg r')\, \dd t'\\
    -\frac{1}{\jm k}\uppartial_\mr{t} \int_S g^\mr{f}(\veg r,\veg r') \uppartial_{\mr t'} j_\mr{t}(\veg r')\, \dd t' \label{eq:TE-EFIE Op}
    \end{multline}
    and the \ac{EFIE} in \eqref{eq:EFIE} becomes
    \begin{equation}
    \TE(j_\mr{t})(\veg r)=-\frac{1}{\eta}e_\mr{t}^\text{inc}(\veg r)\, . \label{eq:TE-EFIE}
    \end{equation}

    \subsection{Low-rank Compression of MoM Matrices}
    
    In the \ac{MoM}, the \acp{SIE} in (\ref{eq:TM-EFIE}) and (\ref{eq:TE-EFIE}) are discretized by approximating $\veg j$ using $N$ (typically localized) basis functions, yielding a system of linear equations $\mat T \vec j = \vec e$.
    Here, $\mat T_{N\times N} \in \{\mat T_{\kern+2pt\text{TM}}, \mat T_{\kern+2pt\text{TE}}\}$ is a dense generalized impedance matrix, $\vec e_{N\times 1} \in \{\vec e_{\text{TM}}, \vec e_{\text{TE}}\}$ contains the tested right-hand side, and $\vec j_{N\times 1}$ is an unknown vector of expansion coefficients, set in accordance with the pertinent equation and operator.

    \begin{figure}
        \centering
        \begin{subcaptionblock}{\linewidth}
        \includegraphics[width=1\linewidth]{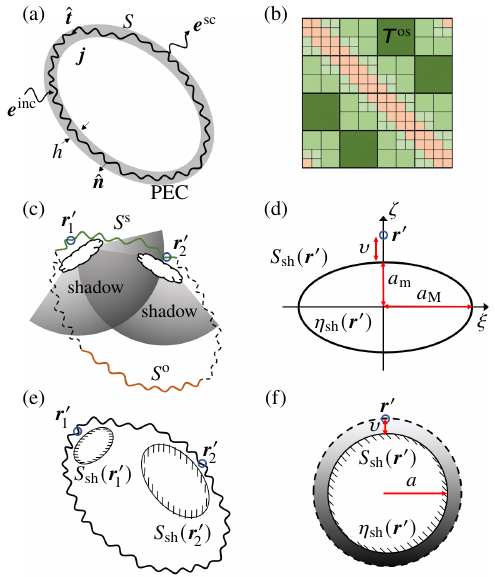}
        \phantomcaption\label{fig:scatteringscenario}
        \phantomcaption\label{fig:hierarchicalblockstructure}
        \phantomcaption\label{fig:attenuatedbroadside}
        \phantomcaption\label{fig:ellipticalshield}
        \phantomcaption\label{fig:elliptialshieldarbit}
        \phantomcaption\label{fig:sourcepointnearcircshield}
        \end{subcaptionblock}
        \caption{\subref{fig:scatteringscenario}  Scattering by an essentially convex surface $S$. \subref{fig:hierarchicalblockstructure} Hierarchical blocks structure of a \ac{MoM} matrix. \subref{fig:attenuatedbroadside} Attenuated broadside radiation in the presence of the shield. \subref{fig:ellipticalshield} Stencil elliptical shield with major and minor axes $a$ and $b$ and a surface impedance  $\eta_\text{sh}$. \subref{fig:elliptialshieldarbit} Two elliptical shields assigned to arbitrary surface source positions $\vec r'_1$ and $\vec r'_2$. \subref{fig:sourcepointnearcircshield} Source point near a circular shield with radius $a$.}
        \label{fig:Fig1}
    \end{figure}

    Algebraic compression-based fast solvers avoid the construction and storage of the entire $\mat T$.
    The matrix is typically partitioned into blocks (as shown in \Cref{fig:hierarchicalblockstructure}), in accordance with the geometric hierarchical clustering of the basis and testing functions \cite{hackbusch1999sparse}.
    A block $\mat T_{N_\mr{o} \times N_\mr{s}}^\text{os}$ of $\mat T$ represents an interaction between clusters of  $N_\mr{s} $ basis functions and  $N_\mr{o} $ testing functions with supports in regions $S^\mr{s}$ and $S^\mr{o}$, respectively. 
    It is deemed “admissible” for compression if $S^\mr{s}$ and $S^\mr{o}$ meet a certain distance and size-based criterion \cite{hackbusch1999sparse}. 
    Commonly, each admissible block $\mat T^\text{os}$ is represented by a \ac{LR} approximation $\mat T^\text{os}\approx \mat A_{N_\mr{o} \times \mathcal{R}} \mat B_{\mathcal{R}\times N_\mr{s}}^\dag$.
    The approximation rank $\mathcal{R}<N_\mr{s}, N_\mr{o}$ is determined by the desired relative error in the representation of $\mat T^\text{os}$.
    For a relative error level $\tau$, the rank $\mc R$ is bounded from below by the number of singular values $\sigma_n$ of $\mat T^\text{os}$ greater than $\sigma_1\tau$.
    This bound is related to the number of degrees of freedom required for accurately representing the underlying continuous interaction between $S^\mr{s}$ and $S^\mr{o}$, as described by the corresponding integral operator \cite{brick2026interpreting, bucci1989degrees}.
    For the conventional operators in \eqref{eq:TM-EFIE} and \eqref{eq:TE-EFIE Op}, a line of sight exist between every pair of points $(\veg r,\veg r')$, which interact as if in free space.
    The line of sight between large subdomains of an essentially convex surface (i.e., a surface defined as a small perturbation from a convex surface) translates to a broadside interaction, for which the rank scale, asymptotically, linearly with $N_\mr{s}$ and $N_\mr{o}$.
    Consequently, the asymptotic scaling with electrical size $kR$ of the computational complexity cannot be reduced .
    
    \subsection{Generalized Source Integral Equations}

    For essentially convex scatterers, defined as a small perturbation in the direction normal to a convex-skull surface (i.e., a large convex shape that can tightly fit within the scatterer, see \Cref{fig:scatteringscenario}), the rank deficiency of moment matrix blocks can be enhanced by modifying the \ac{SIE} kernel in a manner that reduces the effective dimensionality of the interactions.
    For each surface source at point $\veg r'\in S$, an auxiliary contribution to the \ac{SIE} kernel that can be attributed to sources inside $S$ is introduced~\cite{sharshevsky2020direct, zvulun2023generalized, dahan2024fast, dahan2025vector, kalhofer2025fast, kalhofer2026multipole}.
    The auxiliary (dyadic) contribution $\dyd G^\text{aux}(\veg r,\veg r')$  is designed to approximately cancel the field radiated by $\veg j(\veg r')$ into $S$ and, accordingly, to its opposite side (as shown in \Cref{fig:attenuatedbroadside}).
    As a result, broadside interactions are attenuated and end-fire interactions are accentuated, leading to slower scaling of $\mathcal{R}$ with $kR$.
    The generalized (non-physical) sources of the surface distribution $\wveg j(\veg r')$ radiate according to a \ac{MGD}
    \begin{equation}
    \dyd G^\mr{m}(\veg r,\veg r') = \dyd G^\mr{f}(\veg r,\veg r') +\dyd G^\text{aux}(\veg r,\veg r')\, . \label{eq:MGD}
    \end{equation}
    The corresponding \ac{SIE} can be written as 
    \begin{equation}
    \hveg{n}\times\[-\jm k\int_S \dyd G^\mr{m}(\veg r,\veg r')\cdot\wveg j(\veg r')\, \dd t'\] = -\frac{1}{\eta}\n \times \veg e^\text{inc}(\veg r)\, . \label{eq:GSIE}
    \end{equation}
    The \ac{GSIE} counterpart of the TM-EFIE of \eqref{eq:TM-EFIE} can be written, using the scalar $g^\mr{m}(\veg r,\veg r') = g^\mr{f}(\veg r,\veg r')+g^\text{aux}(\veg r,\veg r')$, as
    \begin{equation}
    \TMm(\widetilde j_z)(\veg r) =\TM(\widetilde j_z)(\veg r)+\TMa(\widetilde j_z)(\veg r)=-\frac{1}{\eta}e_\mathit{z}^\text{inc}(\veg r)\, , \label{eq:TM-GSIE}
    \end{equation}
    where 
   \begin{equation}
    \TMa(\widetilde j_z)(\veg r)=-\jm k\int_S g^\text{aux}(\veg r,\veg r') \widetilde j_z(\veg r')\, \dd t' \, . \label{eq:TM-Aux Op}
    \end{equation}
    Similarly, the TE-\ac{GSIE} uses a $2 \times 2$ dyadic $\dyd G^\text{aux}(\veg r,\veg r')$.
    It is described through the vector fields produced by unit currents in each of two primary orthogonal directions $\hveg \xi(\veg r')$ and $\hveg \zeta(\veg r')$ defined at $\vec r'$ (see \Cref{fig:ellipticalshield}), such that 
   \begin{equation}
    \dyd G^\text{aux}(\veg r,\veg r')\cdot\wveg j(\veg r')=\dyd G^\text{aux}\cdot \hveg \xi\widetilde j_\xi(\veg r') + \dyd G^\text{aux}\cdot\hveg \zeta\widetilde j_\zeta(\veg r') \, . \label{eq:TE-Aux}
    \end{equation}
    For the transverse surface distribution $\widetilde j_\mr{t}$, the \ac{GSIE} counterpart of \eqref{eq:TE-EFIE} can be written as
    \begin{equation}
    \TEm(\widetilde j_\mr{t})(\veg r) =\TE(\widetilde j_\mr{t})(\veg r)+\TEa(\widetilde j_\mr{t})(\veg r)=-\frac{1}{\eta}e_s^\text{inc}(\veg r)\, . \label{eq:TE-GSIE}
    \end{equation}
    The \ac{MoM} discretization of these equations results in modified systems $\mat T^\mr{m}\wvec i = \vec v$, where $\mat T^\mr{m} = \mat T + \mat T^\text{aux}$ and $\mat T^\text{aux} \in \lbrace\mat T^\text{aux}_\text{TM},\mat T^\text{aux}_\text{TE} \rbrace$.

    In \cite{dahan2024fast} and \cite{dahan2025vector}, $\dyd G^\text{aux}(\veg r,\veg r')$ corresponds to the fields scattered by a reflective, impenetrable convex shield surface $S_\text{sh}(\veg r')$ when illuminated by sources on $S$.
    The geometric configuration is illustrated in \Cref{fig:ellipticalshield}:
    The surface $S_\text{sh}(\veg r')$ is an ellipse with major and minor axes $a_\mathrm{M}(\veg r')$ and $a_\mathrm{m}(\veg r')$  and has surface impedance $\eta_\text{sh}$; 
    the shield $S_\text{sh}(\veg r')$ is placed such that its minor axis coincides with the direction normal to the convex skull at the point closest to $\veg r'$, at a distance $\upsilon(\veg r')$ to $\veg r'$.
    This placement allows the shielding ellipsoids to lie close to the physical boundary while remaining entirely inside the scatterer.
    In practice, a “stencil auxiliary component” is computed once for a “stencil shield” and then translated and rotated accordingly for each $\veg r'\in S$.
    The  auxiliary contribution can be attributed to induced surface currents on $S_\text{sh}(\veg r')$ or, as in \cite{dahan2024fast} and \cite{dahan2025vector} (and, for completeness, Appendix A) to image sources internal to $S_\text{sh}(\veg r')$.

    For essentially \emph{circular} $S$, effective shielding can be achieved by setting $a_\mathrm{M}=a_\mathrm{m}=a$ (see \Cref{fig:sourcepointnearcircshield}).
    If, in addition, the shields for all $\veg r' \in S$ coincide, the \ac{GSIE} integral operator becomes that of the \ac{GEIE} in \cite{brick2014fastEssentially}.
    The derivation of this early generalized \ac{SIE} assumes a fixed internal scatterer inside $S$ (thus also requiring modification of the \ac{EFIE} right-hand side), as opposed to Love’s equivalence principle.
    For the circular shield case, $\dyd G^\text{m}(\veg r,\veg r')$ can be expressed analytically using cylindrical harmonics. 
    These expressions are utilized in this work for the spectral analysis of the corresponding integral operators.

    \subsection{Conditioning and Calderón Preconditioning Principles}
    \label{subsec:Conditioning and calderon preconditioning principles}
    
    The \acp{GSIE} described by \eqref{eq:TM-GSIE} are, in principle, electric-field equations.
    As such, they may inherit from the \ac{EFIE} various mechanisms that influence their conditioning, with potentially severe implications for the accuracy of the moment solution and for the convergence of iterative solvers \cite{adrian_ElectromagneticIntegralEquations_2021}.
    These include:
    (\textit{i}) Solution non-uniqueness for internal resonance values of $k$. 
    This occurs more frequently with the increase in $k$ for closed objects.
    Conventionally, internal resonances are eliminated by using combined formulations, such as the \ac{CFIE} or \ac{CSIE} \cite{mautz_CombinedsourceSolutionRadiation_1979}.
    (\textit{ii}) Dense-discretization breakdown, that is, the gradual deterioration of the conditioning when the mesh size $h$ decreases \cite{andriulli_MultiplicativeCalderonPreconditioner_2008}.
    (\textit{iii}) High-frequency breakdown, that is, the gradual deterioration of the conditioning when $k$ increases and $hk$ is constant.
    For certain classes of problems, the latter two ill-conditioning mechanisms can be cured by Calderón preconditioning techniques \cite{adrian_ElectromagneticIntegralEquations_2021}.

    The \acp{GSIE} can be viewed as extensions of the \acp{CSIE}~\cite{mautz_CombinedsourceSolutionRadiation_1979}.
    It was suggested (although so far without proof) that their forms in \cite{brick2014fastEssentially,dahan2024fast,dahan2025vector} can even be explicitly designed to represent loss mechanisms that are hypothesized to eliminate internal resonances.
    With regard to the dense-discretization breakdown, we hypothesize, in this work, that the \acp{EFIE} traits inherited by certain \acp{GSIE} can be cured by the established Calderón-type preconditioning techniques.
    For the conventional TM- and TE-\acp{EFIE}, Calderón preconditioners rely on the identity
    \begin{equation}
    \vecop T_\mathrm{\{TE,TM\}}\vecop T_\mathrm{\{TM,TE\}}=-\frac{1}{4}\vecop I + \text{compact}, \label{eq:Claderon-identity}
    \end{equation}
    that is, the $\vecop T_\mathrm{\{TE,TM\}}$ \ac{EFIE} operator preconditions $\vecop T_\mathrm{\{TM,TE\}}$.
    In discretized form, these can be written as 
    \begin{equation}
    {\mat G}^{-\T} \mat T_\mathrm{\{TE,TM\}} {\mat G}^{-1} \mat T_\mathrm{\{TM,TE\}}=-\frac{1}{4}{\mat I} + \text{compact}, \label{eq:Claderon-identity-mat}
    \end{equation}
    where $ {\mat G}$ is a Gram matrix relating the basis and testing functions used for the primal and dual operators.
    For \acp{CFIE} and \acp{CSIE}, the preconditioning operators are typically replaced by their Yukawa kernel counterparts $\vecop Y_\mathrm{\{TE,TM\}}$, defined as the conventional operators with an imaginary wavenumber \cite{contopanagos2002well}, in order to avoid reintroducing interior resonances.
    The exponential decay also allows for sparse storage and fast application of these operators, which aligns with the goals of the \ac{GSIE} approach.
    This suggests that reflective-shield \ac{GSIE}-based operators can be preconditioned with either the $\vecop Y_\mathrm{\{TE,TM\}}$  or $\vecop T_\mathrm{\{TE,TM\}}^{\kern+2pt\mr{m}}$ dual operators.
    To simplify the notation, we denote either choice of dual operator by $\vecop D_\mathrm{\{TE,TM\}}$.
    If, indeed,
    \begin{equation}
    \vecop D_\mathrm{\{TE,TM\}} \vecop T_\mathrm{\{TM,TE\}}^{\kern+2pt\mr{m}}= \\
    -\frac{1}{4}\vecop I + \text{compact}, 
    \label{eq:Claderon-identity_general}
    \end{equation}
    and, in discretized forms,  
    \begin{equation}
    {\mat G}^{-\T} \mat D_\mathrm{\{TE,TM\}} {\mat G}^{-1}\mat T_\mathrm{\{TM,TE\}}^\mr{m}=-\frac{1}{4}{\mat I}+ \text{compact}, \label{eq:Claderon-identity_general-mat}
    \end{equation}
    the resulting formulations are both resonance-free (without the need for a magnetic field equation) and resilient to dense-discretization breakdown.
    For the standard Yukawa-Calderón preconditioned \acp{CFIE} and \acp{CSIE}, it is well known that the high-frequency breakdown persists \cite{boubendir_WellconditionedBoundaryIntegral_2014}.
    The fact that high-frequency breakdown can be cured for limited classes of geometries by using dual operators with locally-modified complex wavenumber kernels \cite{boubendir_WellconditionedBoundaryIntegral_2014} may hint that these local effects could be influenced by the reflective nature of the \ac{GSIE} shields.

    The spectral properties of the \ac{GSIE} and the proposed remedies are investigated both analytically (in closed form, for circular shields and scatterers) and numerically in the following sections. 
    For a matrix $\mat T $, the condition number is defined as $\sigma_1/\sigma_N$. 
    For normal matrices, it is equal to $\max_n\ \abs{\lambda_n} / \min_n\ \abs{\lambda_n} $, where $\lambda_n$ are the \acp{EV} of $\mat T$.
    For the canonical circular geometries used here for the analytical study, the continuous operators admit an explicit eigenvalue decomposition that diagonalizes them in the harmonic basis, with diagonal entries equal to their \acp{EV}.

    \section{Closed-form Spectra for Circularly Symmetric Reflective \acp{GSIE}}\label{sec:Analysis}
        
\begin{figure}[tp]
        \centering
        \includegraphics[width=0.7\linewidth]{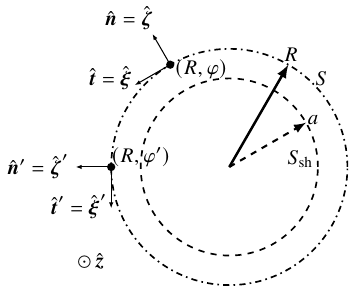}
        \caption{Concentric circular scatterer and shield settings: Dashed line is the shield, dash-dot line the scatterer.
        }
        \label{fig:common problem setting}
    \end{figure}

    The \acp{EV} of the \ac{GSIE} operators are first investigated analytically. 
    This is done for both \ac{TM} and \ac{TE} \acp{SIE} using concentric circular scatterers and coinciding shields, as illustrated in~\Cref{fig:common problem setting}.
    For this configuration, analytical expressions for the \acp{MGF} can be derived. 
    These expressions are then used to derive closed-form expressions for the eigenspectra associated with the harmonic eigenfunctions $\exp(\jm n \varphi)$ $n\in \mathbb{Z}$. 

    \subsection{Modified Green’s Functions}

    Let $S$ and $S_\mathrm{sh}$ be circles of radii $R$ and $a$, respectively.
    On $S_\mathrm{sh}$, a Leontovich~\cite{alexopoulos_accuracy_1975} boundary condition with surface impedance $\eta_\mr{sh} \in \mathbb{C}$ is imposed.
    For this case, the \ac{TM} \ac{MGF} represents the relation between $z$-oriented currents $j_z$ and the $z$-oriented electric field $e_z^\mr{sc}$ they produce in the presence of the internal cylinder. 
    Similarly, the relevant $\xi\xi'$-component of the \ac{TE} \ac{MGD} describes the $\varphi$-directed field produced by a $\hveg \varphi'$-directed current. 
    The derivation of the \ac{MGF} for this case follows a conventional recipe, as that in  \cite[p.~236f.]{harrington_timeharmonic_2001a}, 
    which is outlined for completeness in \Cref{app:modgreentm} of the appendix.

    In polar coordinates, let the source and observer positions be $(\rho',\varphi')$ and $(\rho,\varphi)$.
    For the \ac{TM} formulation, the \ac{MGF} is given by 
    \begin{multline}
        g_\mr{TM}^\mr{m}(\rho, \varphi, \rho', \varphi') =  g^\mr{m}(\rho, \varphi, \rho', \varphi') = \\
        \frac{1}{4\jm}\sum_{n=-\infty}^{\infty} \bigg(\bH_n^{(2)}(k\rho_<)/A^{\mr{TM}}_n + \bJ_n(k\rho_<)\bigg)\bH_n^{(2)}(k\rho_>)\e^{\jm n (\varphi - \varphi')}\,,     \label{eq:TM-GSIE circular shield Green's func}
    \end{multline}
     where
    \begin{equation}
       A^{\mr{TM}}_n  = - \frac{\eta \bH_n^{(2)}(ka) + \jm \eta_\mr{sh} \bH_n^{(2)'}(ka)}{\eta \bJ_n(ka) + \jm \eta_\mr{sh} \bJ_n'(ka)}\,
        \label{eq:An_TM}.
    \end{equation}
    Here, $\rho_{\{>,<\}} = \{\max(\rho,\rho'), \min(\rho,\rho')\}$, $\bH_n^{(2)}(\cdot)$ is the $n$th order Hankel function of the second kind and $\bJ_n(\cdot)$ is the $n$th order Bessel function.
    When both the source and observer lie on $S$, $\rho=\rho'=R$, we have
    \begin{multline}
         g_\mr{TM}^\mr{m}(R, \varphi, R,  \varphi') = \frac{1}{4\jm}\sum_{n=-\infty}^{\infty} \bigg(\bH_n^{(2)}(kR)/A^\mr{TM}_n \\+ \bJ_n(kR)\bigg)\bH_n^{(2)}(kR)\e^{\jm n (\varphi - \varphi')}\,. \label{eq:TM-GSIE circular shield Green's func rho=rho'=R}
    \end{multline}
    For the \ac{TE} formulation, the relevant \ac{MGD} component $g^\mr{m}_{\mr{TE}}=\hveg \varphi \cdot \dyd G^\mr{m} \cdot \hveg \varphi' $ 
    is given by
    \begin{multline}
        g^\mr{m}_{\mr{TE}} (\rho, \varphi, \rho', \varphi')= \frac{1}{4\jm}\sum_{n=-\infty}^{\infty}\bigg( \bH_n^{(2)'}(k\rho_<)/A^{\mr{TE}}_n  \\ + \bJ_n'(k\rho_<)\bigg) \bH_n^{(2)'}(k\rho_>)\e^{\jm n (\varphi - \varphi')}\,,
        \label{eq:TE-GSIE circular shield Green's func}
    \end{multline}
   where
    \begin{equation}
         A^{\mr{TE}}_n  = - \frac{\eta \bH_n^{(2)'}(ka) - \jm \eta_\mr{sh} \bH_n^{(2)}(ka)}{\eta \bJ_n'(ka) - \jm \eta_\mr{sh} \bJ_n(ka)}\,.
           \label{eq:An_TE}
           \end{equation}
     For $\rho=\rho'=R$, this expression reduces to
     \begin{multline}
        g^\mr{m}_{\mr{TE}} (R, \varphi, R, \varphi')= \frac{1}{4\jm}\sum_{n=-\infty}^{\infty}\bigg( \bH_n^{(2)'}(kR)/A^{\mr{TE}}_n  \\ + \bJ_n'(kR)\bigg) \bH_n^{(2)'}(kR)\e^{\jm n (\varphi - \varphi')}\,.
        \label{eq:TE-GSIE circular shield Green's func rho=rho'=R}
    \end{multline}
            
  \subsection{\ac{GSIE} Eigenspectra}
  \label{subsec:Eigenspectra}
 
    For the case of concentric cylinders, the eigenfunction basis for the current $j_z$ consists of the harmonics $\tilde j_z^{(n)}(\varphi) =\e^{\jm n \varphi}$ ($n\in \mathbb{Z}$), which are orthonormal and, therefore, satisfy
    \begin{equation}
        \int_0^{2\pi} \frac{1}{2\uppi}\e^{\jm m (\varphi-\varphi')}\e^{\jm n \varphi'}\, \dd \varphi' = \e^{\jm n \varphi} \delta_{mn}\quad \forall m,n\in \mathbb{Z} \label{eq:Orthogonality}
    \end{equation}
    where $\delta_{mn}$ is the Kronecker delta function.
    Applying the \ac{TM}-\ac{GSIE} operator to any of the eigenfunctions on the circular domain yields
    \begin{equation}
        \TMm(\tilde j_z^{(n)})(\varphi) = -\jm k \int_{0}^{2\uppi} g^\mr{m}_{\mr{TM}}(R, \varphi, R, \varphi')  \tilde j_z^{(n)}(\varphi') R \, \dd \varphi'\,.
    \end{equation}
    Using~\eqref{eq:TM-GSIE circular shield Green's func rho=rho'=R} and~\eqref{eq:Orthogonality}, the \acp{EV} of $ \TMm$ are given by
    \begin{equation}
        \lambda^{\mr{TM}}_n = -\frac{\uppi R k}{2}\(\bH_n^{(2)}(kR)/A^{\mr{TM}}_n +  \bJ_n(kR)\) \bH_n^{(2)}(kR)\,.
        \label{eq:eigenvalues TM-GSIE}
     \end{equation}
    For the \ac{TE}-\ac{GSIE},
    \begin{equation}
        \TEm(\tilde j_\varphi^{(n)})( \varphi) = -\jm k \int_{0}^{2\uppi} g^\mr{m}_{\mr{TE}}(R, \varphi, R, \varphi')  \tilde j^{(n)}_\varphi(\varphi') R \, \dd \varphi'\, .
        \label{eq:TE-GSIE operator reduction}
    \end{equation}
    Substituting~\eqref{eq:TE-GSIE circular shield Green's func rho=rho'=R} into~\eqref{eq:TE-GSIE operator reduction}, the corresponding \acp{EV} are
    \begin{equation}
        \lambda^{\mr{TE}}_n = -\frac{\uppi R k}{2} \(\bH_n^{(2)'}(kR)/A^{\mr{TE}}_n + \bJ_n'(kR)\)\bH_n^{(2)'}(kR) \, .
        \label{eq:eigenvalues TE-GSIE}
    \end{equation}
   
    The expressions in~\eqref{eq:eigenvalues TM-GSIE} and~\eqref{eq:eigenvalues TE-GSIE} provide the spectra of $\TMm$ and $\TEm$. 
    Under Galerkin discretization, the corresponding generalized \acp{EV} of the matrix pairs $(\mat T_\text{TM}^\mr{m},\mat G)$ and $(\mat T_\text{TE}^\mr{m},\mat G)$ converge to these spectra as the mesh size $h \to 0$ and $N \to \infty$.
    Thus, the analytic spectra characterize the \acp{EV} of the Galerkin system matrices and predict their conditioning.
    Next, these spectra are used to examine the sensitivity of these operators to the mechanisms that plague conventional equations and to assess the effectiveness of the preconditioning strategy outlined in \Cref{sec:background}. 

    \subsection{Internal Resonances and Quasi-resonances}
    \label{subsec:InternalResonances}
           
    We examine the hypothesis that the loss mechanisms in the \ac{MGF} serve to eliminate internal resonances. 
    To this end, it should be determined whether $\lambda^{\mr{TM}}_n = 0$ and $\lambda^{\mr{TE}}_n = 0$ are possible \acp{EV}. For the \ac{TM} operator, this translates to
    \begin{equation}
        -\frac{\uppi R k}{2}\(\bH_n^{(2)}(kR)/A^{\mr{TM}}_n +  \bJ_n(kR)\) \bH_n^{(2)}(kR) = 0\, , 
        \label{eq:lambda_v^M = 0}
    \end{equation}
    implying that both the real and imaginary parts of the term in the parentheses should be zero.
    Expressing $\eta_\mr{sh} = \eta_\mr{R} + \jm \eta_\mr{I}$ using $\eta_\mr{R}, \eta_\mr{I} \in \R$ and denoting $x=ka$, the condition in \eqref{eq:lambda_v^M = 0} can be written as 
    \begin{align}
        C_\Re &:= \eta_\mr{R}\( \bJ_n(Rx/a)\bY_n'(x) - \bJ_n'(x)\bY_n(Rx/a) \) = 0\label{eq: C1}\,,\\
        C_\Im &:= \frac{\eta_\mr{I}}{\eta} \(\bJ_n(Rx/a)\bY_n'(x) - \bJ_n'(x)\bY_n(Rx/a)\) \notag \\
        &\quad - \( \bJ_n(Rx/a)\bY_n(x) -  \bJ_n(x)\bY_n(Rx/a)\) = 0\label{eq: C2}\,,
    \end{align}
    These equations can be viewed as a homogeneous linear system for the unknowns $\bJ_n(\alpha x)$ and $\bY_n(\alpha x)$, where $\alpha = R/a$, which takes the following matrix form
    \begin{equation}
       \begin{pmatrix}
           \eta_\mr{R}\bY_n'(x) & -\eta_\mr{R}\bJ_n'(x) \\
           \frac{\eta_\mr{I}}{\eta} \bY'_n(x) - \bY_n(x) & -\frac{\eta_\mr{I}}{\eta} \bJ_n'(x) + \bJ_n(x)
       \end{pmatrix}
       \begin{pmatrix}
           \bJ_n(\alpha x) \\
           \bY_n(\alpha x)
       \end{pmatrix}
       = \begin{pmatrix}
           0\\0
       \end{pmatrix}\,.
    \end{equation}
    This homogeneous system can only have a nontrivial solution if the determinant of the left-hand side matrix  is zero.
    However, for $\eta_\mr{R}\neq 0$, the determinant is proportional to the Wronksian
    \begin{equation}
        W\{\bJ_n,\bY_n\}(x) = \bJ_n(x)\bY_n'(x) - \bJ_n'(x)\bY_n(x) =  \frac{2}{\pi x}\,,\label{eq:WronskianLHS}
    \end{equation}
    which has no nulls.
    A similar condition can be derived for the \ac{TE} problem and its \acp{EV}. 
    This suggests that these \acp{GSIE} are, in principle, resonance-free for non-zero $\eta_{\mr R}$.
    
\begin{figure}
        \centering
        \begin{subcaptionblock}{\linewidth}
            \includegraphics[]{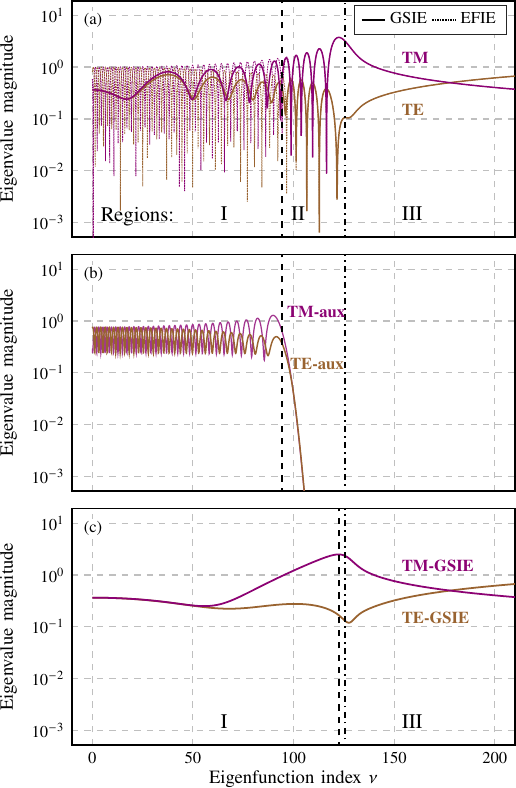}
    
            \phantomcaption\label{fig:GSIE-lam0.5-q5}
            \phantomcaption\label{fig:GSIEgaux-gf-lam0.5-q5}
            \phantomcaption\label{fig:GSIE-lam0.5-q0.5}            
        \end{subcaptionblock}

        \caption{Concentric circular scatterer--shield configuration: Closed-form \ac{EV} spectra of \ac{TM} and \ac{TE} operators, with spectral region identification. (a) \ac{GSIE} spectra for $\upsilon= 5\lambda$, (b) \ac{GSIE} spectra for $\upsilon =0.5\lambda$, and (c) Auxiliary and free space contributions \ac{GSIE} spectra for $\upsilon= 5\lambda$.}
        \label{fig:GSIE spectra}
    \end{figure}

    The non-zero~\ac{GSIE} \acp{EV} may still be arbitrarily close to zero.
    For example, consider the \ac{EV} spectra for the choice $R=10\lambda$ ($\lambda = \qty{1}{\m}$) and $a = R - \upsilon$ and two choices of  $\upsilon \in \mathbb{R}^+$, shown in~\Cref{fig:GSIE spectra}.
    The number of indices displayed corresponds to the expected value of $N$ with $h \approx\lambda /10$. 
    To visualize the underlying spectral trends, the closed-form eigenvalue expressions are evaluated for continuous orders $\nu\in\R^+$ rather than only for integer indices $n\in\mathbb{Z}$.
    Since cylindrical functions of negative integer order can be expressed in terms of those of positive order~\cite[Sec.~3.53]{watson1944treatise}, only non-negative indices are considered in the following.

    When examining the spectrum, it is convenient to partition it into three regions: $\nu<ka$ (Region~I), $ka<\nu<kR$ (Region~II), and $\nu>kR$ (Region~III). 
    The spectrum may exhibit fundamentally different behavior in each of these regions.      
    For example, in the case of $\upsilon=5\lambda$ (\Cref{fig:GSIE-lam0.5-q5}), combs of sharp dips are observed for the \acp{GSIE} only in Region~II, for both the \ac{TM}- and \ac{TE}-\ac{GSIE}.
    These dips, which are not exactly zero, are attributed to radial “quasi-resonances”.
    For comparison, the \ac{EFIE} exhibits these dips, which are proper resonances here, in both Regions~I and~II.

    The behavior in Region~II can be explained by examining the contribution of the auxiliary Green’s function to the spectrum (see~\Cref{fig:GSIEgaux-gf-lam0.5-q5}).
    It exhibits nulls at different $\nu$ values than those of the free-space contribution in Region~II.
    Consequently, when added to the free-space contribution, it cancels its internal resonances.
    For $\nu>ka$, the auxiliary contribution enters its static/elliptic regime, where it decays rapidly, even if not to zero. 
    For larger $\upsilon$, this leaves Region~II governed by the resulting quasi-resonant behavior of the free-space kernel.
    Since the free-space contribution is already smooth in Region~III, only Region~II is affected by this transition.
    For $\upsilon=0.5\lambda$ (\Cref{fig:GSIE-lam0.5-q0.5}), the width of Region~II decreases and does not allow for the existence of the quasi-resonance dips.
    Generally, the quasi-resonances can be avoided by setting $\upsilon < \lambda $ (which conveniently aligns with the desire for deep shadow that translates to enhanced compressibility).
 
    Although sharp quasi-resonance dips occur only in Region~II, large $\upsilon$ values produce bounded oscillations in Region~I.
    The local minima in this region can be associated with lossy modal behavior that does not exist for lossless choices of $\eta_\mr{sh}$ (for which the distinction between Regions~I and~II disappears).

    \subsection{Dense-Discretization Breakdown and Proposed Cure}
    \Cref{fig:GSIEgaux-gf-lam0.5-q5} also shows the dense-discretization breakdown of the conventional \ac{TM} and \ac{TE}~\acp{EFIE}. 
    The \ac{TM}-\ac{EFIE} \acp{EV} cluster at zero with rate $\mc O(h)$ (Region~III may be regarded as extending to infinity to capture the behavior when $h \to 0 $), whereas the \ac{TE}-\ac{EFIE} \acp{EV} approach infinity with rate $\mc O (1/h)$~\cite{adrian_ElectromagneticIntegralEquations_2021}. 
    The operators $\TMa$ and $\TEa$ are compact perturbations of $\TM$ and $\TE$, which means that their \acp{EV} cluster at zero for $n \to \infty$. 
    Consequently, the auxiliary contribution in Region~III (see~\Cref{fig:GSIEgaux-gf-lam0.5-q5}) is negligible and the dense-discretization properties of the \ac{EFIE} are inherited by the \ac{GSIE}. 
    This implies that the Calderón-type preconditioning techniques commonly used for \acp{EFIE} can be used with the \acp{GSIE} to combat the dense-discretization breakdown.
    As a result, the \ac{EV} magnitudes cluster around $1/4$ in Region~II, as observed in~\cref{fig:hf spectra}.
     
    \subsection{High-Frequency Behavior}
    The conventional \ac{CSIE} also suffers from the growth of the condition number with increasing frequency while keeping $kh$ constant, commonly referred to as the “high-frequency breakdown”.
    For the \ac{TM} and \ac{TE} problems, it is attributed to the local maxima and minima near $\nu=ka$ in \Cref{fig:GSIE spectra} (also termed the glancing spectrum \cite{antoine_ImprovedSurfaceRadiation_2006}). 
    These local extrema scale faster with $k$ than the remaining \acp{EV}, leading to a deterioration of the condition number at a rate $\propto k^{1/3}$.
    This behavior is not cured by the multiplicative Yukawa-Calderón preconditioning.

    \begin{figure}
        \centering
        \begin{subcaptionblock}{\linewidth}
            \includegraphics[]{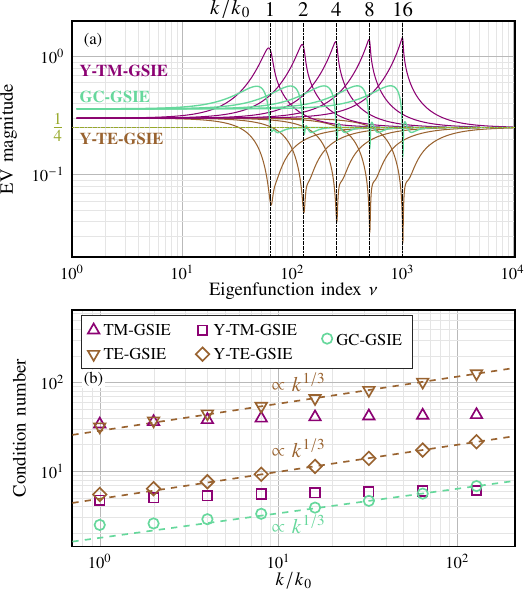}

            \phantomcaption\label{fig:hf spectra}
            \phantomcaption\label{fig:hf breakdown}
            
        \end{subcaptionblock}
        \caption{Concentric circular scatterer--shield configuration: (a) Spectra of the preconditioned operators with increasing frequency (b) Condition numbers for the different \ac{GSIE} formulations with increasing frequency.}
        \label{fig:hf behavior}
    \end{figure}
        
    For the \acp{GSIE}, the scaling of the extrema with $k$ is demonstrated in~\Cref{fig:hf breakdown}.
    The \ac{TM}-GSIE operator shows only a mild growth of the condition number, with a trend of $\propto k^{0.05}$.
    This indicates that the shields influence the local wave mechanism responsible for the high-frequency breakdown. 
    When the Yukawa \ac{TE}-\ac{EFIE} operator is applied as a preconditioner to the \ac{TM}-\ac{GSIE}, the resulting equation (Y-\ac{TM}-\ac{GSIE}) is expected to inherit this favorable behavior, while not exhibiting a dense-discretization breakdown.
    This is thanks to the  exponentially decaying, highly-localized Yukawa kernel, which does not result in a decreasing dip in $k$ for the corresponding  \ac{TE}-\ac{EFIE} operator, as also demonstrated in~\Cref{fig:hf breakdown}.
    However, the \ac{TE}-\ac{GSIE} operator shows a growth rate approaching $k^{1/3}$, similar to conventional operators.
    Accordingly, neither a Yukawa-Calderón preconditioned nor a \ac{TM}-\ac{GSIE} Calderón preconditioned~\ac{TE}-\ac{GSIE} is free from high-frequency breakdown.
    Although both approaches improve the conditioning (Y-\ac{TE}-\ac{GSIE} and \acs{GC}-\ac{GSIE}, respectively), the $k^{1/3}$ growth persists.
    Likewise, the \ac{TE}-\ac{GSIE} operator as a Calderón preconditioner for the \ac{TM}-GSIE is expected to introduce a worse scaling trend.
    Indeed, with both preconditioning types, the $k^{1/3}$ growth rate is observed, though the effectiveness of the GSIE-based Calderón preconditioner for the \ac{TE}-\ac{GSIE} is greater. 
    As a preconditioner, the \ac{TE}-\ac{GSIE} operator is slightly superior to the Yukawa-Calderón alternative until a cross-over is reached.
    
    These phenomena will be examined numerically, including non-circular scatterers and stencil shield-based \acp{GSIE}, in the next section.      

    
    \section{Numerical Results}\label{sec:NumResults}
    The spectral properties of \ac{GSIE} \ac{MoM} matrices and the effectiveness of preconditioning are examined numerically.
    The first examples confirm that the spectra of the moment matrices converge to the \acp{EV} of the continuous operators for concentric circular scatterer--shield configurations. 
    A second set of examples verifies that spectral properties of the circular-shield \acp{GSIE} are similar to those based on snugly fitting convex stencil shields \cite{dahan2024fast,dahan2025vector}.  
    Then, for circular scatterers and elliptical shields, the dependence of the conditioning on the mesh size and the effectiveness of Calderón-type preconditioners for the \ac{GSIE} are examined. 
    Lastly, the dependence of the conditioning on the frequency is studied.
    Throughout this section, unless otherwise specified, elliptical shields with an aspect ratio $a_\text{M}=1.5a_\text{m}$ and an impedance $\eta_\text{sh}=(1+\jm)\qty{1.37e2}{\ohm}$ are used. 
    The \ac{MoM} discretization uses polygonal meshes with piecewise-linear (hat) basis and testing functions, defined for the mesh nodes. 
    For this choice, the entries of $ \mat G$  are 
    \begin{equation}
        \sbr{\mat G}_{ij}=\int \Lambda_i \Lambda_j \,\dd s = \begin{cases} 
        \frac{1}{6}h & i=j+1\ \text{or}\ i=j-1 \\
        \frac{2}{3}h & i=j \\
        0 & \text{else}
        \end{cases}, \label{eq: Gram-matrix}
    \end{equation}
    A four-point Gaussian quadrature rule is used.
    Singular integrals are evaluated analytically via small-argument approximation of their integrands.
    
    \subsection{Spectra of the Continuous and Discretized Circular Shield GSIE}
    \begin{figure}
    \centering
    \begin{subcaptionblock}{\columnwidth}
        \includegraphics[width=\columnwidth]{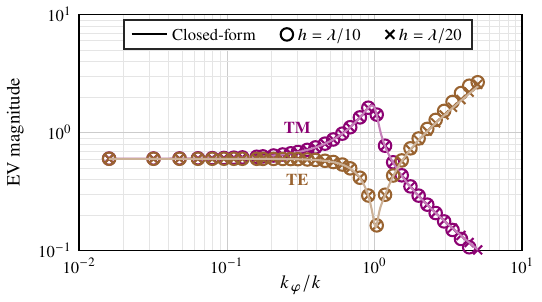}
    \end{subcaptionblock}
    \caption{Concentric circular scatterer--shield configuration: Closed-form and \ac{MoM} matrix spectra of the \ac{TM}- and \ac{TE}-\ac{GSIE}}
    \label{fig:Analytical_Numerical_Circle_GSIE}
    \end{figure}
    For concentric circular scatterer--shield configurations, the closed-form spectra in \eqref{eq:eigenvalues TM-GSIE} and \eqref{eq:eigenvalues TE-GSIE} are presented as functions of $k_\varphi/k = n/(kR)$ and compared to the spectra of the corresponding \ac{MoM} matrices after left multiplication by $\mat G^{-1}$.
    The results, shown in \Cref{fig:Analytical_Numerical_Circle_GSIE}, 
    are for $R=10\lambda$, $\upsilon=0.3\lambda$,
    and two mesh sizes $h \in \{\lambda/10,\lambda/20 \}$.
    For the continuous operators, which have infinitely many \acp{EV} with known asymptotic behavior, only the first $N/2$ positive-index eigenvalues are plotted.
    The \ac{MoM} spectra agree well with the closed-form spectra in Regions~I and~II.
    In Region~III, they converge to the closed-form \acp{EV} as the discretization is refined, with good agreement observed already for $h=\lambda/20$.
    For this choice of $\upsilon=0.3\lambda<\lambda$ and $\eta_\mr{sh}$, the spectra do not exhibit dips and are free of quasi-resonances. 

    \subsection{Stencil Shield-Based \ac{GSIE} Spectra}

    \begin{figure*}
    \centering
    \begin{subcaptionblock}{\linewidth}
        \includegraphics[width=\linewidth]{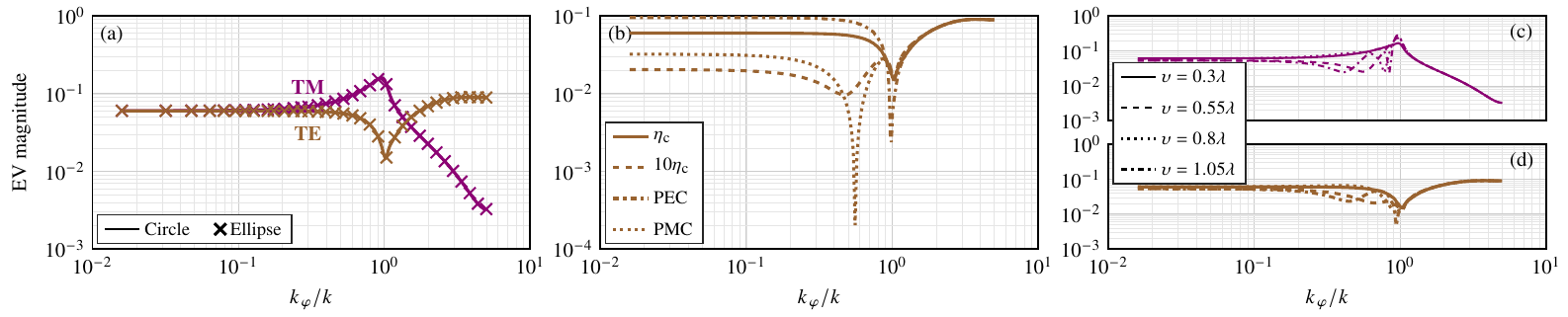}
        \phantomcaption\label{fig:GSIE_Circle_Ellipse}
        \phantomcaption\label{fig:GSIE_Ellipse_varImpedance}
        \phantomcaption\label{fig:TM_GSIE_varDistance}
        \phantomcaption\label{fig:TE_GSIE_varDistance}
    \end{subcaptionblock}
    \caption{(a) \ac{GSIE} moment matrix spectra for circular and ellipse shields. (b) Comparison of spectra for various $\eta_\text{sh}$ values. (c,d) \ac{GSIE} spectra for several values of $ \upsilon/\lambda$ for the (c) \ac{TM} and (d) \ac{TE} formulations.}
    \label{fig:CircleEllipseSpectrum}
    \end{figure*}

    Next, we examine whether the spectral characteristics of the \acp{GSIE} moment matrices with concentric circular shields are maintained with snugly fitting elliptical shields.
    This is done for the circular scatterer of \Cref{fig:common problem setting}, shields with proportions $a_\text{M}=1.5a_\text{m} = 0.6R$ and $\upsilon = 0.3\lambda$.
    \Cref{fig:GSIE_Circle_Ellipse} shows the spectra for the default $\eta_\text{sh}$ value.
    Excellent agreement is observed for both \ac{TM} and \ac{TE} formulations.
    Even though they do not correspond to a \ac{GSIE} with a single well-defined internal scatterer for describing the auxiliary contribution, the elliptical-shield-based \acp{GSIE} maintain the practically resonance-free nature of their circular-shield counterparts.
    In particular, the agreement of the extrema near $k_\varphi/k=1$ implies that the local wave mechanisms that govern them, which are influenced by the curvature of $S$, $\upsilon$, and $\eta_\mr{sh}$, are accurately mimicked by the elliptical-shield-based \acs{GSIE}.

    To examine the influence of $\eta_\mr{sh}$, the spectra for various values of the shield impedance are compared for the \ac{TE}-\ac{GSIE} formulation with $\nu=0.3\lambda$.
    The comparison includes the values $\eta_\mr{sh} \in \{\eta_\text{c},10\eta_\text{c}\}$ where $\eta_\text{c}=(1+\jm)\qty{1.37e2}{\ohm}$ as well as \ac{PEC} and \ac{PMC} boundaries. 
    The results are shown in \Cref{fig:GSIE_Ellipse_varImpedance}. 
    For the \ac{PEC} and \ac{PMC} boundaries, the spectra exhibit deeper and sharper minima, approaching nulls, within $k_\varphi/k < 1 $. 
    These nulls correspond to resonances and occur regardless of the sub-wavelength value of $\upsilon$. 
    Notably, for the \ac{PMC} case, the resonance appears deep within what was classified as Region~I, suggesting the merging of Regions~I and~II. 
    For the lossy boundaries, Region~II is small and no high-Q quasi-resonances are observed near $k_\varphi/k = 1$. 
    The impedance does influence the existence of broad modal minima. 
    For the choice $\eta_\mr{sh} =(1+\jm)\qty{1.37e3}{\ohm}$, such a mode appears well within Region~I. 
    Because that broad minimum is deeper than the local minimum near $k_\varphi/k=1$, it governs the condition number until higher frequencies are reached.   
    In \Crefrange{fig:TM_GSIE_varDistance}{fig:TE_GSIE_varDistance}, the influence of $\upsilon$ is examined for the \ac{TM}- and \ac{TE}-\acp{GSIE}. 
    The spectra for $\upsilon/\lambda \in \{0.3, 0.55, 0.8, 1.05\}$ are presented. 
    Consistent with the results in \Cref{subsec:InternalResonances}, 
    as $\upsilon$ is increased, additional modal minima appear in Region~I.
    Once $\upsilon=\lambda$ is exceeded, a radial quasi-resonance appears for the \ac{TE} case.
 
    \subsection{GSIE Dense-Discretization Breakdown and Cures}
    The influence of mesh refinement on the spectrum and conditioning of the moment matrix and the effectiveness of Calderón-type preconditioning is examined next.

    \begin{figure*}
    \centering
    \begin{subcaptionblock}{\linewidth}
        \includegraphics[width=\linewidth]{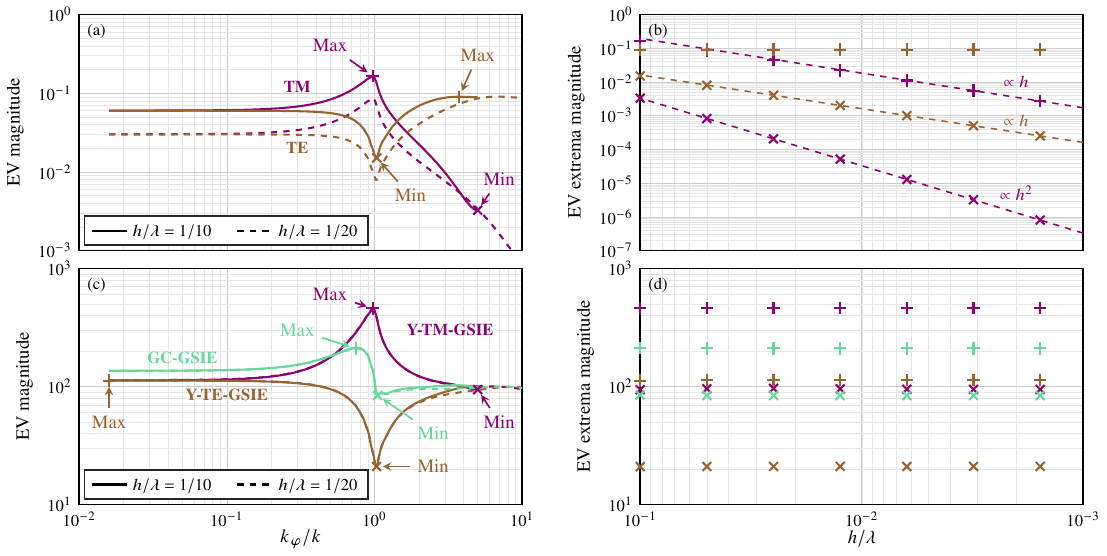}
        \phantomcaption\label{fig:DenseDisc_Spectra_Ellipse}
        \phantomcaption\label{fig:DenseDisc_Spectra_Extrema}
        \phantomcaption\label{fig:DenseDisc_Spectra_Preconditioned}
        \phantomcaption\label{fig:DenseDisc_Spectra_PreconditionedExtrema}
    \end{subcaptionblock}
    \caption{(a) Spectra of elliptical shield \ac{TM}- and \ac{TE}-\acp{GSIE}. (b) Extrema of the spectra in (a) as functions of $h$. (c)  Spectra of \ac{TM}- and \ac{TE}-\acp{GSIE} preconditioned using the Yukawa and GSIE kernel duals. (d) Extrema of the spectra in (c) as functions of $h$.}
    \label{fig:Dense_discretization_Spectra}
    \end{figure*}
    
    For the circular scatterer and the elliptical shields of \Cref{fig:GSIE_Circle_Ellipse}, 
    \Cref{fig:DenseDisc_Spectra_Ellipse} presents the \ac{TM}- and \ac{TE}-\ac{GSIE} moment-matrix spectra for $h/\lambda= 1/10$ and  $h/\lambda= 1/20$. 
    The spectra scale with $h$ for $k_\varphi$-values in Regions~I and~II, and approach fixed asymptotic values for large $k_\varphi$. 
    The extrema are indicated by markers on the $h/\lambda= 1/10$ plots. 
    The scaling of the condition number with $h$ is determined by the ratio of these extrema.
    \Cref{fig:DenseDisc_Spectra_Extrema} shows the scaling of these extrema with $h/\lambda$.
    For each of the formulations, their ratios scale asymptotically as $\propto 1/h$, suggesting that both \acp{GSIE} inherit the dense-discretization breakdown of their conventional \ac{EFIE} counterparts.
    
    Following \eqref{eq:Claderon-identity_general} and \eqref{eq:Claderon-identity_general-mat}, the effectiveness of the low-rank-compressible \ac{GSIE} dual and its Yukawa-kernel counterpart in curing the dense-discretization breakdown is examined next. 
    \Cref{fig:DenseDisc_Spectra_Preconditioned} shows the spectra of the preconditioned operators for $h/\lambda= 1/10$ and $h/\lambda=1/20$.
    The extrema for the two values of $h$ align, showing no deterioration of the condition number.
    This behavior is shown to be maintained for broader $h/\lambda$-value ranges in \Cref{fig:DenseDisc_Spectra_PreconditionedExtrema}, where the extrema indicated in \Cref{fig:DenseDisc_Spectra_Preconditioned} keep constant levels in the entire $h/\lambda$-range, suggesting that the proposed remedies are, indeed, effective.

    \begin{figure}
        \centering
        \begin{subcaptionblock}{\linewidth}
            \includegraphics[width=\columnwidth]{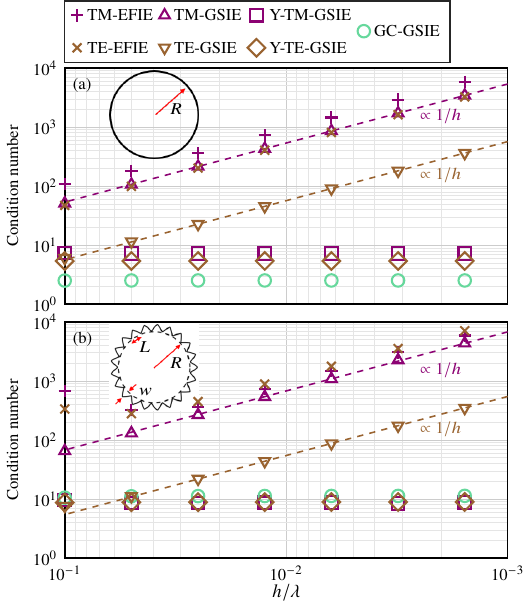}
            \phantomcaption\label{fig:DenseDisc_Condition_NumberCircular}
            \phantomcaption\label{fig:DenseDisc_Condition_NumberCorrugate}
        \end{subcaptionblock}
        \caption{Condition numbers of the moment matrices as functions of $h$ for nonpreconditioned and preconditioned operators, for (a) circular $S$ and (b)  corrugated circle $S$. A dashed line indicates a circular skull.}
        \label{fig:Dense discretization condition number}
    \end{figure}

    The condition numbers of the preconditioned and nonpreconditioned moment matrices that were analyzed in \Cref{fig:Dense_discretization_Spectra} are shown in 
    \Cref{fig:DenseDisc_Condition_NumberCircular}.
    The preconditioned \ac{GSIE} operators maintain constant conditioning with $h/\lambda$.
    Notably, using the \ac{GSIE}-dual results in the lowest condition number among all studied cases. 
    The study is repeated for the corrugated circular cylinder in the inset of~\Cref{fig:DenseDisc_Condition_NumberCorrugate}.
    The jagged corrugation has a period of $L=1.5\lambda$ and height $w=0.6\lambda$ above a circular skull surface. 
    Here, $\upsilon\in[0.2, 0.8]\lambda$, with the shields are placed at a fixed distance of $0.2\lambda$ from the circular skull of $S$ (see \Cref{fig:DenseDisc_Condition_NumberCorrugate} for a visualization of the configuration).
    The condition numbers in~\Cref{fig:DenseDisc_Condition_NumberCorrugate} show a similar behavior to that in~\Cref{fig:DenseDisc_Condition_NumberCircular} for the smooth circular $S$.
  
    \subsection{GSIE High-Frequency Breakdown}
    Lastly, the influence of the frequency on the conditioning of the \ac{GSIE} \ac{MoM} matrices is studied.
    This is done for a circular scatterer of $R=10\lambda_0$ (where the wavelength $\lambda_0$ corresponds to a wavenumber $k_0$).
    For each frequency, the value $\upsilon=0.3\lambda$ is kept to avoid internal quasi-resonances.

   \begin{figure*}
    \centering
    \begin{subcaptionblock}{\linewidth}
        \includegraphics[width=\textwidth]{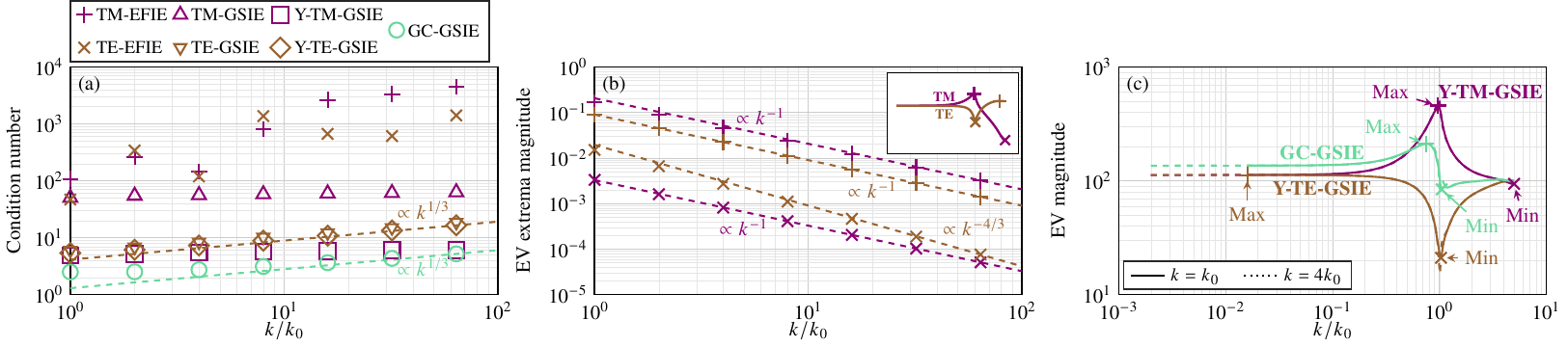}
        \phantomcaption\label{fig:HighFrequencyCond}
        \phantomcaption\label{fig:HighFrequencyExtremaMag}
        \phantomcaption\label{fig:HighFrequencyMag}
    \end{subcaptionblock}
    \caption{(a) Moment matrix condition numbers scaling with $k$ for various formulations. (b) Extrema of the \ac{TM}- and \ac{TE}-\ac{GSIE} spectra as functions of $k$. (c) Spectra of \ac{TM}- and \ac{TE}-\ac{GSIE} preconditioned using Yukawa and \ac{GSIE} kernel duals for two $k$ values.}
    \label{fig:High frequency condition number & spectrm}
    \end{figure*}
    
    \Cref{fig:HighFrequencyCond} shows the moment matrix condition number scaling with $k$ for the formulations considered in \Cref{fig:DenseDisc_Condition_NumberCorrugate}.
    Note that, for the \acp{EFIE}, the displayed values are calculated by excluding the minimal \ac{EV}, which can be associated with an internal resonance.
    The \acp{EFIE} exhibit the expected deterioration of the conditioning with $k$, as does the nonpreconditioned \ac{TE}-\ac{GSIE}, with a trend that approaches $k^{1/3}$.
    However, the nonpreconditioned \ac{TM}-\ac{GSIE} condition number exhibits a slow $\propto k^{1/20}$ scaling.
    When applying the Yukawa-kernel preconditioner to the \ac{TM}-\ac{GSIE}, the condition number is lowered and deterioration with~$k$ remains practically invisible.
    Alternatively, when using the \ac{TE}-\ac{GSIE} dual as a preconditioner, the $k^{1/3}$ scaling remains despite the significant reduction in the condition number resulting from the elimination of the dense-discretization breakdown.
    Similar behavior is observed when applying the Yukawa dual to the \ac{TE}-\ac{GSIE}.
    This is in agreement with that demonstrated in \Cref{fig:hf breakdown} for the concentric circular scatterer--shield configuration. 
    It is due to the faster decrease with $k$ of the \ac{TE}-\ac{GSIE} glancing spectrum minimum (see~\Cref{fig:HighFrequencyExtremaMag}), unlike for the \ac{TM}-\ac{GSIE}, where the extrema scale similarly.
    Hence, while the \ac{TM}-\ac{GSIE} behaves well following the application of the Yukawa-dual preconditioner, the \ac{TE}-\ac{GSIE} is not similarly cured of its high-frequency breakdown. 
    Furthermore, if the \ac{TE}-\ac{GSIE} is used to precondition the \ac{TM}-\ac{GSIE}, while lowering the condition number for moderate frequencies, it ruins its inherent relative resilience to the high-frequency breakdown, causing an eventual $k^{1/3}$ scaling of the condition number.
    For the selected shield parameters, these behaviors of the preconditioned operators can be seen already in the transition from $k=k_0$ to $k=4k_0$ in~\Cref{fig:HighFrequencyMag}. 
    It is worth noting that in some cases, for example when $\eta_\mr{sh}=10\eta_\text{c}$ as in \Cref{fig:CircleEllipseSpectrum} for the \ac{TE}-\ac{GSIE}, a broad local minimum that scales similarly to the maximum may govern the conditioning for sufficiently low frequencies, making it appear indifferent to $k$, until the glancing spectrum minimum becomes the dominant one and the high-frequency deterioration appears.

    \section{Conclusion}\label{sec:Conclusion}
    This work establishes the spectral foundations of reflective \acp{GSIE}.
    Although the auxiliary kernels fundamentally alter the radiation pattern in order to improve the compressibility of moment matrix blocks, they remain compact perturbations of the corresponding \ac{EFIE} operators.
    Consequently, while the \acp{GSIE} inherit the dense-discretization breakdown of their \ac{EFIE} counterparts, they respond well to the same Calderón-type preconditioning principles that cure it for \acp{EFIE}.
    The analysis also reveals that the auxiliary kernels influence the mechanisms governing the high-frequency spectrum.
    For the \ac{TM} formulation, the reflective shield fundamentally changes the dominant spectral branch, resulting in an unusually mild growth of the condition number that is preserved under Yukawa-Calderón preconditioning.
    The \ac{TE} formulation, in contrast, retains the conventional high-frequency deterioration.
    These observations demonstrate that auxiliary kernels affect not only the compressibility but also the spectral mechanisms governing iterative solution. 
    
    Perhaps most importantly, the theoretical predictions derived for the continuous circular operators remain valid for discretized operators and practical elliptical stencil shield-based \acp{GSIE}.
    The resulting \ac{GSIE} formulations therefore combine three properties that have previously been studied largely independently: freedom from internal resonances, resilience to dense-discretization breakdown through Calderón-type preconditioning, and preservation of the low-rank structure required by modern fast solvers.
    Harnessing these properties for the treatment of more complex and 3-D objects, specifically via geometrically more-enabling \ac{GSIE} kernel designs (e.g., \cite{kalhofer2025fast, kalhofer2026multipole}), is an avenue for further research. 
	

	\appendices
    \section{Reflective Shields MGF and MGD}
    In the implementation chosen for this work, the auxiliary contributions to the Green's functions are produced by sets of $M$ auxiliary sources with complex weights $M_m$. 
    Since the shield geometry is identical for every source location up to translation, the weights $\{M_m\}_{m=1}^M$ need only be computed once for an arbitrary reference point $\veg r_0$ and a stencil shield $S_\mr{sh}(\veg r_0)$, within which the auxiliary sources are placed at points $\veg r_m$. 
    The auxiliary sources radiate according to the free-space Green's kernel.
    For example, for the \ac{TM}-\ac{GSIE}, the auxiliary contribution is produced by $z$-directed current lines and the \ac{MGF} is given by
    \begin{equation}
          g^\mr{m}(\veg r,\veg r_0) = g^\mr{f}(\veg r,\veg r_0)+\sum_{m=1}^M M_m g^\mr{f}(\veg r,\veg r_m)\, .
        \label{eq:TMg E-field}
    \end{equation}
    The weights should be set such that the impedance boundary condition is satisfied at the observation points $\veg r'' \in S_\mr{sh}(\veg r_0)$.   
    Denoting
    \begin{equation}
        \veg g_{H}^{\mr{f}} (\veg r, \veg r') = 
        \nabla_\mr{t}\times \(\hveg z g^\mr{f}(\veg r,\veg r')\)
        \label{eq:TM gH-field}
    \end{equation}
    and $\hveg t''_\text{sh} (\veg r')$ the unit vector tangent to $S_\text{sh}(\veg r')$ at $\veg r''$, this amounts to requiring that
        \begin{multline}
        \ -\jm k \eta \[g^\mr{f}(\veg r'',\veg r_0)+\sum_{m=1}^M M_m g^\mr{f}(\veg r'',\veg r_m)\]\ = \\
        \eta_\text{sh} \hat{\veg t}''_\text{sh} (\veg r_0) \cdot \[\veg g _{H}^{\mr{f}} (\veg r'', \veg r_0) + \sum_{m=1}^M M_m \veg g_{H}^{\mr{f}} (\veg r'', \veg r_m)\]\, , \\ \forall \veg r'' \in S_\mr{sh}(\veg r_0)\, , 
        \label{eq:TMg impedance condition}
    \end{multline}
    This requirement can be imposed at $M$ or more points on $S_\mr{sh}(\veg r_0)$, leading to a (potentially, overdetermined) system of equations that can be solved in the least-squared sense. 

    The \ac{TE}-\ac{GSIE} uses a $2\times2$ dyadic $\dyd G^\text{aux}(\veg r,\veg r')$ for two orthogonal components of the surface current (see~\eqref{eq:TE-Aux}).
    The $ {\{\xi, \zeta\}}$ contribution is produced by a corresponding set of auxiliary sources with weights $ M^{\{\xi, \zeta\}}_m$.   
    Here, following \cite{dahan2025vector}, magnetic current line sources are used. 
    Denoting 
    \begin{equation}
        \veg g_{E}^{\mr{f}} (\veg r,\veg r') = + \frac{1}{k^2}\nabla_\mr{t}\times \(\hveg z g^\mr{f}(\veg r,\veg r')\)
        \label{eq:magneticSourceEfield}\, ,
    \end{equation}
    the stencil auxiliary component to the \ac{MGD} by a $\hveg \beta \in \{ \hveg \xi, \hveg \zeta \}$ directed unit current is written as
    \begin{equation}
        \dyd G^\text{aux}(\veg r,\veg r_0)\cdot\hveg \beta  = \sum_{m=1}^M M^{\beta }_m \veg g_{E}^{\mr{f}} (\veg r,\veg r_m)\, ,
        \label{eq:TE-Aux decomposition}
    \end{equation}
    implying that $ M^{ \{\xi, \zeta\} }_m$ are measured in \si{\per\meter}. 
    The $z$-directed magnetic stencil field is produced by a $\hveg \beta \in \{ \hveg \xi, \hveg \zeta \}$ oriented unit current source at $\veg r_0$. 
    Defining
    \begin{equation}
         g_{H\beta}^\mr{f}(\veg r,\veg r') = \hveg z \cdot [\nabla g^\mr{f}\times \hveg \beta]
        \label{eq:currentSourceHfield}
    \end{equation}
    and
    \begin{equation}
        g_{H\beta}^\text{aux}(\veg r,\veg r') = \sum_{m=1}^M M^\beta_m  g^{\mr{f}} (\veg r,\veg r_m)\, ,
        \label{eq:magneticSourceHfield}
    \end{equation}
    the impedance boundary condition implies that
    \begin{multline}
       -\jm k\eta \hveg t''_\text{sh} (\veg r_0) \cdot \[ \dyd G^\mr{f}(\veg r'',\veg r_0) \cdot \hveg \beta+\sum_{m=1}^M M^\beta_m \veg g_{E}^{\mr{f}} (\veg r'',\veg r_m)\] = \\
        -\eta_\text{sh}\[g_{H\beta}^\mr{f}(\veg r'',\veg r_0) +\sum_{m=1}^M M^\beta_m g^{\mr{f}} (\veg r'',\veg r_m)\]\, , \quad \veg r'' \in S_\mr{sh}(\veg r_0)\, .
        \label{eq:TEg impedance condition}
    \end{multline}

	\section{\ac{GSIE} \acp{MGF}  for Circular Shields}\label{app:modgreentm}
    The derivation of the \acp{MGF} follows~\cite[Sec.~(6.5.2)]{jin_theory_2015} closely. 
    For simple sources placed, without loss of generality, at $\varphi'=0$, the auxiliary contribution is the electric field scattered by the impedance boundary circular shield \cite{alexopoulos_accuracy_1975}.

    In this case, the Leontovich boundary condition~\eqref{eq:TMg impedance condition} reduces to
    \begin{equation}
        -\jm k \eta \( g^\mr{f}  + g^\mr{aux} \) = \eta_\mr{sh} \hveg{\varphi} \cdot \( \veg g_H^\mr{f} + \veg g_H^\mr{aux} \) \,.
        \label{eq:app_Leontovich_TM}
    \end{equation}
    Each of the components can be expressed as an infinite harmonic series such that
    \begin{align}
         - \jm k \eta \gf &= -\frac{\omega \mu}{4}\sum_{n=-\infty} ^{\infty} \bJ_n(k\rho)\bH_n^{(2)}(k\rho')\e^{\jm n\varphi}\, , \label{eq:app_Ez_inc_jin}\\
  - \jm k \eta \gaux &=  -\frac{\omega \mu}{4} \sum_{n=-\infty} ^{\infty} d_n \bH_n^{(2)}(k\rho) \bH_n^{(2)}(k\rho') \e^{\jm n\varphi}, \label{eq:app_Ez_sc_jin}\\
         \hveg{\varphi} \cdot \gfH &= \frac{\jm k}{4}\sum_{n=-\infty}^{\infty} \bJ_n'(k\rho) \bH_n^{(2)}(k\rho') \e^{\jm n \varphi}\, ,
         \label{eq:a5}\\
     \hveg{\varphi} \cdot \gauxH &= \frac{\jm k}{4}\sum_{n=-\infty} ^{\infty} d_n \bH_n^{(2)'}(k\rho) \bH_n^{(2)}(k\rho') \e^{\jm n\varphi}\, ,
        \label{eq:app_Hz_sc_jin}
    \end{align}
    Enforcing the the boundary condition in \eqref{eq:app_Leontovich_TM} at $\rho = a$ (see~\Cref{fig:common problem setting}), $d_n = 1 /A_n^\mr{TM}$ and
     $A_n^\text{TM}$ are the coefficients in~\eqref{eq:An_TM}.
    The form in \eqref{eq:TM-GSIE circular shield Green's func} is obtained by summing the free-space and auxiliary components above. 

    For the \ac{TE} problem, the boundary condition in \eqref{eq:app_Leontovich_TM} reduces to
    \begin{equation}
      -\jm k\eta \hveg \varphi \cdot \[ \dyd G^\mr{f} \cdot \hveg \varphi'+\dyd G^\mr{aux} \cdot \hveg \varphi'\] = 
        -\eta_\text{sh}\[g_{H\varphi}^\mr{f} + g_{H\varphi'}^\text{aux}\]\, .
        \label{eq:bpp_impedance_bc}
    \end{equation}
    Here, only the $\varphi$ field component produced by a $\varphi'$-directed current source is of interest.
    Without loss of generality, the source is assumed to be placed at $\varphi'=0$ and $y$-directed.
    The corresponding field components are expressed as
    \begin{align}
        -\jm k\eta \hveg \varphi \cdot \dyd G^\mr{f} \cdot \hveg \varphi' &= -\frac{\omega \mu}{4} \sum_{n=-\infty}^{\infty}  \bJ_n'(k\rho) \bH_n^{(2)'} (k\rho') \e^{\jm n \varphi}\, ,\label{eq:bpp_Ephi_inc}\\
        -\jm k\eta \hveg \varphi \cdot \dyd G^\mr{aux} \cdot \hveg \varphi' &= -\frac{\omega \mu}{4}\sum_{n=-\infty}^{\infty} d_n \bH_n^{(2)'}(k\rho) \bH_n^{(2)'}(k\rho') \e^{\jm n \varphi}\, ,\\
        g_{H\varphi}^\mr{f} &= \frac{\jm k}{4}\sum_{n=-\infty}^{\infty}
        \bJ_n(k\rho) \bH_n^{(2)'}(k\rho') \,\e^{\jm n \varphi}\, ,\\
         g_{H\varphi}^\text{aux} &=  \frac{\jm k}{4}\sum_{n=-\infty}^{\infty} d_n \bH_n^{(2)}(k\rho) \bH_n^{(2)'}(k\rho') \e^{\jm n \varphi}\, ,
    \end{align}
    Substituting these into~\eqref{eq:bpp_impedance_bc} yields $d_n = 1/A_n^\mr{TE}$, with the values $A_n^\mr{TE}$ in~\eqref{eq:An_TE}.
    The expression in \eqref{eq:TE-GSIE circular shield Green's func} is obtained by summation of the electric field contributions.



%

	\ifCLASSOPTIONcaptionsoff
	  \newpage
	\fi

	

	\printbibliography

\end{document}